\documentclass{article}
\usepackage[french,UKenglish]{babel} 
\usepackage[T1]{fontenc} 
\usepackage[utf8]{inputenc}
\usepackage{fullpage} 
\newcommand\systF{\mathbf{F}}
\usepackage{amssymb} 
\usepackage{parallel} 
\newcommand\imp{\mathbin{\rightarrow}}
\newcommand\paragraphit[1]{\paragraph{\bf\itshape #1}\it} 

\title{\it On the system $\systF$ as a glue language\\
for natural-language compositional-semantics\\[1em]
\rm Du système $\systF$ comme langage d'assemblage\\ 
pour la sémantique compositionnelle du langage naturel} 
\author{\textsc{Christian RETOR\'E} (Université de Bordeaux, LaBRI \& INRIA)}

\begin{document}
\maketitle
\center 
\begin{tabular}{p{0.50\textwidth}p{0.47\textwidth}}
\paragraph{Résumé}
Afin de modéliser dans un cadre compositionnel des phénomènes de pragmatique lexicale comme ceux étudiés par
\cite{cruse1986lexical,Pus95} 
et plus particulièrement 
par Nicholas Asher dans 
\cite{asher-webofwords}, 
notre équipe a utilisé dans divers travaux 
\cite{BMRjolli,Mery2011,MPR2011taln,MPR2011cid} le système 
$\systF$ de Jean-Yves Girard (1971) 
\cite{Girard71,GLT88} pour construire les formules logiques représentant le sens des énoncés 
--- là où d'autres auteurs comme 
\cite{Luo2011lacl} utilisent plutôt la théorie des types (1984) de Per Martin-Löf. \cite{ML84} 
Nous expliquons dans cette note les raisons de notre préférence.  
&
\paragraphit{Abstract} 
In order to model in compositional framework some phenomena of lexical pragmatics as the ones studied by 
\cite{cruse1986lexical,Pus95} 
and especially by Nicholas Asher in 
\cite{asher-webofwords}, 
several contributions developed in our team \cite{BMRjolli,Mery2011,MPR2011taln,MPR2011cid}  did use the 
system 
$\systF$ of Jean-Yves Girard (1971) 
\cite{Girard71,GLT88} 
to construct logical formulae expressing the meaning
---  while others like \cite{Luo2011lacl} prefer to use Per Martin-Löf's type theory \cite{ML84}. 
In this note we explain the motivations supporting our preference. 
\end{tabular} 

\begin{tabular}{p{0.50\textwidth}p{0.47\textwidth}}
\paragraph{Qu’est-ce que le système $\systF$?}  Introduit dans \cite{Girard71}, 
(une bonne référence est \cite{GLT88} — en anglais et disponible sur Internet — une  autre est \cite{Girard2006pa1}), 
ce système est la logique propositionnelle intuitionniste du second ordre. On y quantifie sur les propositions: on prouve 
des formules du genre 
$\forall p\  (p\imp p)\imp (p \imp p)$ noté $\Pi p.\ (p\imp p)\imp (p \imp p)$ 
et les preuves s'écrivent comme 
des termes avec un $\Lambda p$ terme représentant la règle de quantification sur les types (introduction du $\Pi\ p$ ci-dessus), par exemple: 
$\Lambda p.\ \lambda f^{p\imp p}\lambda x^p  f(f(x))$ qui est du type $\Pi p.\ (p\imp p)\imp (p \imp p)$ — la variable de type sur laquelle on quantifie ne doit pas être libre dans le type d'une variable libre. Symétriquement, un terme générique $u$ de type $\Pi p.\ T[p]$ peut être spécialisé sur un type $A$: le terme correspondant s'écrit $u\{A\}$, et il est de type $T[A]$. Par exemple si $u$ trie de liste de $p$-objets pour tout $p$, $u\{A\}$ trie des listes de $A$-objets. 

La beta réduction usuelle sur $\lambda$ est
$(\lambda x^T.\ u) t^T\rightsquigarrow u[x:=t]$ mais  celle sur $\Lambda$ fait pareil: $(\Lambda \alpha.\ u)\{T\} \rightsquigarrow u[\alpha:=T]$. 
&
\paragraphit{What is system $\systF$?}  Introduced in \cite{Girard71},  (a good reference freely available on the Web is \cite{GLT88} 
and another is \cite{Girard2006pa1}), this system is intuitionistic second order propositional logic. In this system one can quantify over propositions and can prove formula like 
$\forall p\  (p\imp p)\imp (p \imp p)$ written $\Pi p.\ (p\imp p)\imp (p \imp p)$ 
and proofs are written as terms with $\Lambda p$ denoting the rule introducing quantification over types  --- the type variable over which one quantifies should not be free in the type of a free variable. Symmetrically, a generic term $u$ de type $\Pi p.\ T[p]$ can be specialised to a type $A$: the corresponding term, which is of type $T[A]$ is denoted by $u\{A\}$.   For instance, if $u$ sorts lists of $p$-objects for all $p$, $u\{A\}$ sorts lists of $A$-objects. 

Usual beta reduction concerns $\lambda$: 
$(\lambda x^T.\ u) t^T\rightsquigarrow u[x:=t]$ but beta reduction on $\Lambda$ works just the same: $(\Lambda \alpha.\ u)\{T\} \rightsquigarrow u[\alpha:=T]$. 
\end{tabular} 

\begin{tabular}{p{0.50\textwidth}p{0.47\textwidth}}
\paragraph{Capacité expressive} 
On remarque facilement que la quantification universelle sur les types permet de définir toutes sortes d'opérateurs de constructions de types, 
notamment les plus courants, en particulier:
\begin{itemize} 
\item le produit cartésien de deux types $A\times B\equiv \Pi X. (A\imp (B\imp X)\imp X$
\item la quantification  existentielle sur les types $\exists p A[p] \equiv \Pi q.\ (\Pi p.\ (A[p]\imp q))\imp q$
\item les booléens $\Pi X. X \imp X\imp X$
\item les entiers $\Pi X. (X \imp X)\imp (X\imp X)$
\item les listes d'objets de type $\alpha$:  $\Pi X. X \imp (\alpha \imp X\imp X) \imp X$
\item ... d'où son succès pour ML, CaML etc. 
\end{itemize} 
On peut représenter toutes les fonctions récursives totales dont on peut prouver la totalité dans l'arithmétique de Peano du second ordre 
(qui pour ce genre d'énoncés coïncide avec l'arithmétique de Heyting). 
&
\paragraphit{Expressive power} 
Universal quantification over types let you define lots of type constructors, 
including the common ones, in particular:
\begin{itemize} 
\item cartesian product $A\times B\equiv \Pi X. (A\imp (B\imp X)\imp X$
\item existential quantification over types $\exists p A[p] \equiv \Pi q.\ (\Pi p.\ (A[p]\imp q))\imp q$
\item booleans $\Pi X. X \imp X\imp X$
\item integers $\Pi X. (X \imp X)\imp (X\imp X)$
\item lists of objects of type $\alpha$:  $\Pi X. X \imp (\alpha \imp X\imp X) \imp X$
\item ... hence its success with ML, CaML etc. 
\end{itemize} 
A way to characterise is expressive power is to say that in this system one can define exactly all the total recursive functions 
whose totality can be proved in second order Peano arithmetic (classical or intuitionistic it makes no difference). 
\end{tabular} 

\begin{tabular}{p{0.50\textwidth}p{0.47\textwidth}}
\paragraph{La logique d’assemblage n'est pas la logique où formuler la sémantique} 
Le système $\systF$ (logique propositionnelle du second ordre) sert de logique d'assemblage, mais il assemble les formules que l'ont veut, et 
ce peut bien sûr être des formules de la logique du premier ordre --- même en utilisant des quantifications du second ordre, cela n'a pas de rapport. 
Un résultat dans la thèse de Bruno Mery \cite{Mery2011} est le suivant: 
si les constantes du $\lambda$-calcul sont celles de la logique du premier ordre 
(resp. d'ordre $n$, d'ordre $\omega$), alors les $\lambda$-termes de $\systF$ 
normaux de type $t$ sont les formules de la logique du premier ordre  (resp. d'ordre $n$, d'ordre $\omega$). Cela est bien sûr à rapprocher de la sémantique de Montague usuelle où le $\lambda$-calcul simplement typé (logique propositionnelle intuitionniste) avec deux types de base $e$ et $t$ 
permet d'exprimer et d'assembler les formules de la logique  du premier ordre, d'ordre $n$  ou d'ordre $\omega$.
&
\paragraphit{The glue logic has no relation with the logic used for formulating semantics} 
The system $\systF$ (intuitionistic \emph{second} order propositional logic)  is used as a glue language, but it glues whichever formulae one likes. 
and it can be \emph{first} order formulae --- eventhough one uses second order quantifications, these are unrelated languages. 
A result in the PhD of Bruno Mery \cite{Mery2011} garantees  that when $\lambda$-calculus constants 
are the ones of first (resp. $n$, $\omega$order logic, the normal $\lambda$-terms of $\systF$ of type $t$ are formulae of first (resp. $n$, $\omega$) order logic.
This is  close to usual Montague semantics 
where simply typed $\lambda$-calculus (intuitionistic propositional logic) 
with two base types $e$ and $t$ is used to express and to glue formulae of first order logic, $n$ order logic and even $\omega$ order logic. 
\end{tabular} 

\begin{tabular}{p{0.50\textwidth}p{0.47\textwidth}}
\paragraph{Sous-typage}
Comme j'en ai parfois fait la remarque, le sous-typage est plus ou moins incompatible le système $\systF$ malgré diverses tentatives comme \cite{Cardelli94Fsubtyping}. 
Mais je ne suis pas sûr que le sous-typage au sens de la programmation fonctionnelle soit ce dont on a besoin pour la sémantique 
dans le style de \cite{Pus95,asher-webofwords}. 
Un calcul typé avec sous-typage permet, à partir d'inclusions sur les types atomiques, 
de dériver systématiquement les inclusions sur les types complexes. 
On peut se demander si, sémantiquement,  les inclusions pertinentes de types composés $(a\imp b)$ sont toutes celles issues 
de celles sur $a$ et de celles sur $b$. 
Je ne suis pas sur que le sous typage des verbes, par exemple, se déduise du sous-typage 
sur leurs arguments, sujet, objet, etc.   
La classification des aliments ou des "mangeurs" et des "aliments" induit-elle une classification 
des variantes de manger (\emph{bouffer, dévorer, déguster}) ? Je ne crois pas. 
On veut donc plutôt veut un mécanisme de sous-typage indépendant
sur chaque famille de types, et ce genre de "sous typage" n'a rien à voir avec une relation de sous-typage au sens de la programmation fonctionnelle.  On peut même remarquer que la langue ne permet pas toutes les inclusions ontologiques, et qu'il n'est même pas clair que les inclusions arbitraires du langage  définissent  un ordre. 
&
\paragraphit{Subtyping} 
As I sometimes said, sub-typing is more or less incompatible with system $\systF$ despite some attempts like \cite{Cardelli94Fsubtyping}
But I am not sure that subtyping in its functional-programming meaning is what we are looking for when dealing with the semantic questions of 
 \cite{Pus95,asher-webofwords}. A typed calculus with subtyping automatically 
 derives inclusions between complex types from inclusions between base types. 
 We should wonder whether it is semantically sound that inclusions between complex types like $a\imp b$ 
 are all the ones derived from inclusions on $a$ and inclusions on $b$. 
 For instance, I am not sure that subtyping on verb types derive from their arguments, subject, object, etc. 
Does classifications of  "food"and  "eaters" provide a classification of "eating" verbs (\emph{swallow, taste, appreciate})?
I do not think so. We are rather looking for an independent "subtyping" mechanism for each type (or type family) and this kind of subtyping has nothing to do with the standard notion subtyping at work in functional programming. One can even observe that language does not allow all the ontological inclusions and that it is not so clear that the idiosyncratic linguistic inclusions define an order. 
\end{tabular} 

\begin{tabular}{p{0.50\textwidth}p{0.47\textwidth}}
\paragraph{Intérêt de la quantification sur les types} 
C'est assez pratique. En particulier il n'y a qu'un quantificateur universel, qu'on spécialise, et non un par sorte. 
Sans le $\Lambda$ il faudrait $\forall_a: (a\imp t)\imp t$ pour chaque type $T$, y compris pour les types complexes. 
Ici, un seul suffit: $\forall: \Pi a.\ ((a\imp t)\imp t)$, on l'applique ensuite à $T$ pour obtenir le quantificateur sur $T$: 
 $\forall\{T\}:((T\imp t)\imp t)$. (On eput même avoir une constante  $\forall: \Pi a.\ t$, sans  \textsl{type raising}). 
De même on peut avoir un "il exists", un "la plus part" etc. qui ensuite sont spécialisés pour chaque type particulier. 
On peut aussi traiter dans ce cadre le $\epsilon$ ou $\tau$ de Hilbert qui sélectionne un individu satisfaisant une propriété dès qu'il en existe 
— notamment utilisé pour modéliser le voyageur fictif  
$\tau: \Pi a.\  ((a\imp t)\imp e)$. \cite{MPR2011taln,MPR2011cid} 

On utilise des transformations et des instanciations sur des types composés, 
comme par exemple  pour commuer un chemin en voyageur qui l'emprunte, en événement etc. 
Par exemple, comme ci-dessus, on peut opèrer sur le type raising, $(voie \imp t) \imp  t$ qui devient $((humain \imp t) \imp t$.  \cite{MPR2011taln,MPR2011cid} 
&
\paragraphit{Quantifying over types is great} 
At least, it is useful. In particular we have a single universal quantifier, that is specialised, and not one per sort. 
Without the $\Lambda$ one would need $\forall_a: (a\imp t)\imp t$ for every type $T$, including the complex types. 
Here, one is enough: $\forall: \Pi a.\ ((a\imp t)\imp t)$, then it is applied to a type $T$ to obtain the quantifier over $T$ objects: 
$\forall\{T\}:((T\imp t)\imp t)$. 
(One can even have a constant $\forall: \Pi a.\ t$, without type raising). 
Similarly one can have a "there exists", a "most of" that are specialised to every particular type. 
The very same trick can be used for Hilbert's $\epsilon$ and $\tau$ selecting an object satisfying a property as soon as there is one
--- we already used it, for instance, in the virtual traveller 
$\tau: \Pi a.\  ((a\imp t)\imp e)$. \cite{MPR2011taln,MPR2011cid}

We use transformations and instantiations on complex types, 
for instance to turn a path into a traveller following it, into an event, etc. 
For instance, as above, we can operate on the type raising: 
 $(path \imp t) \imp  t$ qui devient $((human \imp t) \imp t$. 
\end{tabular}

\begin{tabular}{p{0.50\textwidth}p{0.47\textwidth}}
\paragraph{Questions de complexité}
En syntaxe avoir une classe de grammaires la plus petite possible permet d'avoir des algorithmes d'analyse plus performants.
En sémantique, il n'en est rien: que doit faire au minimum la logique d'assemblage? on n'en sait rien...
mais si on le savait cela améliorerait-il la complexité algorithmique de l'analyse sémantique? Non, pas du tout! Ce sont plutôt les choix des différents "sens" de chaque mot qui rendent le processus complexe. 
Le processus d'analyse est l'analyse syntaxique  que vous préférez (et là il y a des questions de complexité), puis 
on construit un $\lambda$-terme avec les $\lambda$-termes sémantiques du lexique (un choix dont la  complexité ne saurait être améliorée) et les règles syntaxiques, et on le réduit, 
et bien évidemment on n'a pas de beta réduction problématique
--- par exemple nécessitant un nombre exponentiel d'étapes de beta réduction en fonction de la longueur du terme à beta réduire. 
&
\paragraphit{Complexity issues} 
In syntax, having a class of grammars as restricted as possible, allows more efficient parsing algorithms. 
In semantics, it  does not make sense: what is the minimum that the glue logic should do? No one knows, but even if this was known, would it help to reduce the algorithmic complexity of semantic analysis? Not at all! 
The complexity is rather due to the choice, for each word, of its relevant meaning among the possibilities 
provided by the lexicon. The process of analysis is simply syntactic analysis, and the syntactic analysis of a phrase can be the one you prefer (here is a complexity issue), then a $\lambda$-term is constructed following the syntax 
with the chosen semantic $\lambda$-terms  from the lexicon (a choice whose complexity cannot be improved), then it is reduced and obviously the reduction is unproblematic --- one does not get a term whose number of reduction steps is exponential with respect to the length of the term. 
\end{tabular} 

\begin{tabular}{p{0.50\textwidth}p{0.47\textwidth}}
\paragraph{Comparatif entre la théorie des types de Martin-Löf et le système $\systF$ pour la sémantique}
Certains,  comme Zhaohui Luo \cite{Luo2011lacl},  utilisent plus volontiers la théorie des types de Martin-Löf (TT) dont il existe plusieurs variantes, l'une des premières étant \cite{ML84}. 
Une discussion avec un expert des deux systèmes, Thierry Coquand (qui a étendu $\systF_\omega$ aux types dépendants, pour faire le calcul des constructions \cite{constructions88}
et qui depuis près de vingt ans travaille sur TT) a achevé de me convaincre: comme logique d'assemblage 
il vaut mieux choisir le 
système $\systF$ pour sa simplicité formelle (4 règles) et conceptuelle (une seule sorte de type, un seul genre de jugements), 
car on n'a aucun problème de complexité 
dans la réduction --- comme dit plus haut, on ne va pas chercher les fonctions à croissance ultra rapide du système $\systF$ pour assembler des formules. 

Si on prend une version faible de TT, celle-ci  est, du point de vue des fondements --- reconnaissons-le --- plus simple. 
La normalisation des termes du système, permet d'établir sans trop de peine la cohérence du système
(s'il existait une preuve du faux il en existerait une qui soit normale, et on peut voir assez facilement qu'il n'y en a pas).
Pour TT il faut le schéma de compréhension (qui permet de définir un ensemble à partir d'une formule $\{X | P(X)\}$  
pour les formule $Pi_1^1$ c'est-à-dire pour les formules de la forme 
"Pour tout ensemble $X$ Il existe un ensemble $Y$..." tandis que pour $\systF$ le schéma de compréhension pour toutes les formules. 
&
\paragraphit{Comparing Martin-Löf's type theory and system $\systF$ for semantics} 
Some, like Zhaohui Luo \cite{Luo2011lacl} rather us Martin-Löf's type theory (TT),
a system with several variants, and early one being in \cite{ML84}. 
A converging discussion with an expert in both, namely Thierry Coquand (who extended $\systF_\omega$ to dependent types, yielding the calculus of constructions \cite{constructions88}, working with TT for twenty years) 
achieved to convince me: as a glue logic, system $\systF$ should be preferred, because of its formal complexity, 
because of its conceptual simplicity (a single sort of types and a single sort of de judgements), because we have no complexity in reduction --- as earlier said, we are not going to use extremely fast growing functions and the like when glueing formulae. 

If one uses a weak variant of TT, this system is admittedly simpler as far as foundations are concerned. 
Indeed, coherence is derived from normalisation (if there were a proof of a false statement, there would be a normal proof, but it is easily seen that there cannot be a normal proof of something false). 
For a weak variant of TT, the comprehension axiom (which defines a set $\{X | P(X)\}$ from a formula $P$)  is needed for every $\Pi_1^1$ formula $P$ (formulae like "for every set $X$ there  exists a set $Y$....."). 
For system $\systF$ the comprehension axiom is needed for every formula. 
\end{tabular} 

\begin{tabular}{p{0.50\textwidth}p{0.47\textwidth}}
\paragraph{\textsl{Records} et types dépendants}
Quant au type de données populaire \textsl{records}, il s'agit simplement d'un produit (définissable dans  dont on nomme les projections, un peu comme les attributs d'une base de données. 
Ils peuvent aussi se voir comme un cas particulier d'une famille d'ensembles indexée par un autre ensemble.
Ce sont alors donc des types dépendants particuliers, 
et certains pensent que, d'une manière générale,  les types dépendants 
sont utiles à la sémantique compositionnelle, ce qui, personnellement, ne me saute pas aux yeux. Néanmoins, pour les inconditionnels  des types dépendants,  
ou si quelqu'un en montrait l'utilité,  on peut les ajouter au système $\systF$.
On peut même les ajouter à $F_\omega$ ce qui donne un système bien connu: le calcul des constructions de Coquand-Huet,
implémenté dans l'assistant de preuve Coq. \cite{constructions88,BC04} 
& 
\paragraphit{Records and dependent types} 
Regarding the popular data type known as "records", they are just products with named projections,
a bit like attributes in a data base. 
They can be viewed as particular indexed family of sets.
Hence they are particular dependent types,
and some thinks that dependent types are relevant for compositional semantics,
also personnally I am not convinced. 
In any case, if one is keen on dependent types, or if one show they are useful at semantics, 
they could be added to system $\systF$. They can even be added to 
$F_\omega$ 
yielding a famous system: the calculus of constructions of Coquand-Huet, implemented in the Coq proof assistant. \cite{constructions88,BC04}
\end{tabular} 

\small
\bibliographystyle{plain}

\end{document}